\documentclass{gtart}

\def\ifplaintex{\expandafter\ifx\csname documentclass\endcsname\relax}

\ifplaintex 
\else
\fi

\def\gtm{{\mathsurround=0pt\it $\cal G\mskip-2mu$eometry \&\ 
$\cal T\!\!$opology $\cal M\mskip-1mu$onographs}}    

\def\gtp{{\mathsurround=0pt\it $\cal G\mskip-2mu$eometry \&\ 
$\cal T\!\!$opology $\cal P\!$ublications}}  

\def\recd{{\small Received:\qua\receiveddate\ifx\reviseddate\relax
\else\qquad Revised:\qua\reviseddate\fi\par}} 

\def\volumenumber#1{\def\thevolumenumber{#1}}
\def\volumeyear#1{\def\thevolumeyear{#1}}
\def\volumename#1{\def\thevolumename{#1}}
\def\papernumber#1{\def\thepapernumber{#1}}
\def\pagenumbers#1#2{\def\startpage{#1}\def\finishpage{#2}}
\def\published#1{\def\publishdate{#1}}
\def\received#1{\def\receiveddate{#1}}
\def\revised#1{\def\reviseddate{#1}}
\def\accepted#1{\def\accepteddate{#1}}

\let\\\par
\let\thevolumenumber\relax\let\thepapernumber\relax
\let\thevolumeyear\relax\let\startpage\relax
\let\finishpage\relax\let\publishdate\relax\let\receiveddate\relax
\let\reviseddate\relax\let\accepteddate\relax\let\theasciititle\relax
\let\theasciiauthors\relax
\let\theasciiabstract\relax

\let\theerratum\relax\let\theasciiemail\relax
\let\theshortauthors\relax\let\theshorttitle\relax

\def\startpage{1}\def\finishpage{15}\def\thepapernumber{77}

\volumenumber{2}
\volumename{Proceedings of the Kirbyfest}
\volumeyear{1999}

\long\def\maketitlep{   

\count0=\startpage

\gtm\nl        
{\small Volume \thevolumenumber: \thevolumename\nl 
\ifx\theerratum\relax\else Erratum \erratumnumber\nl\fi
Pages \startpage--\finishpage\nl}

\vglue 0.1truein   

{\parskip=0pt\leftskip 0pt plus 1fil\def\\{\par\smallskip}{\ifplaintex\large
\else\Large\fi\bf\thetitle}\par\medskip}   
\vglue 0.05truein 

{\parskip=0pt\leftskip 0pt plus 1fil\def\\{\par}{\sc\theauthors}
\par\medskip}%
 
\vglue 0.03truein 

{\small\leftskip 25pt\rightskip 25pt{\bf Abstract}\stdspace\theabstract

{\bf AMS Classification}\stdspace\theprimaryclass
\ifx\thesecondaryclass\relax\else; \thesecondaryclass\fi\par
{\bf Keywords}\stdspace \thekeywords\par}\vglue 7pt

}   

\font\phead=cmsl9 scaled 950
\font=cmsl9 scaled 1050
\font\pnum=cmbx10 scaled 913
\font\lnum=cmbx10 
\font\pfoot=cmsl9 scaled 950
\font=cmsl9 scaled 1050
\ifplaintex
\headline{\vbox to 0pt{\vskip -4.5mm\line{\small\phead\ifnum
\count0=\startpage ISSN 1464-8997 (on line)
1464-8989 (printed) \hfill {\pnum\folio}\else\ifodd\count0\def\\{ }%
\ifx\theshorttitle\relax\thetitle\else\theshorttitle\fi\hfill{\pnum\folio}
\else\def\\{ and }{\pnum\folio}\hfill\ifx\theshortauthors\relax\theauthors
\else\theshortauthors\fi\fi\fi}\vss}}
\footline{\vbox to 0pt{\vglue 0mm\line{\small\pfoot\ifnum\count0=\startpage
Published \publishdate:\qua\copyright\ \gtp\hfill\else
\gtm, Volume \thevolumenumber\ (\thevolumeyear)\hfill\fi}\vss
}}
\else
\makeatletter
\def\@oddhead{{\small\lhead\ifnum\count0=\startpage ISSN 1464-8997 (on line)
1464-8989 (printed) \hfill {\lnum\number\count0}\else\ifodd\count0
\def\\{ }\ifx\theshorttitle\relax \thetitle \else\theshorttitle\fi\hfill
{\lnum\number\count0}\else\def\\{ and }{\lnum\number\count0}
\hfill\ifx\theshortauthors\relax 
\theauthors\else\theshortauthors\fi\fi\fi}}\def\@evenhead{@oddhead}
\def\@oddfoot{\small\lfoot\ifnum\count0=\startpage Published \publishdate:\qua\copyright\ \gtp\hfill\else
\gtm, Volume \thevolumenumber\ (\thevolumeyear)\hfill\fi}
\def\@evenfoot{@oddfoot}
\makeatother
\fi

\let\maketitlepage\maketitlep
\let\makeshorttitle\maketitlepage
\let\maketitle\maketitlepage

\newwrite\gtoutfile
\long\gdef\makeheadfile{  
{\def\\{, }\def\s{ }
\immediate\openout\gtoutfile head.xxx
\immediate\write\gtoutfile{To: math@arxiv.org}
\immediate\write\gtoutfile{Subject: put OR rep NNNNN:ppppp}
\immediate\write\gtoutfile{--text follows this line--}
\immediate\write\gtoutfile{Proxy-for: \ifx\theasciiauthors\relax
\theauthors\else\theasciiauthors\fi\s<\ifx\theasciiemail\relax\theemail\else\theasciiemail\fi>}
\immediate\write\gtoutfile{\noexpand\\}
\immediate\write\gtoutfile{Authors: \ifx\theasciiauthors\relax
\theauthors\else\theasciiauthors\fi}
{\def\\{ }\immediate\write\gtoutfile{Title: \ifx\theasciititle\relax
\thetitle\else\theasciititle\fi}}
\immediate\write\gtoutfile{Subj-class: GT or SG, GR etc}
\immediate\write\gtoutfile{MSC-class: \theprimaryclass\ifx\thesecondaryclass\relax\else, \thesecondaryclass\fi}
\immediate\write\gtoutfile{Journal-ref: Geom. Topol. Monogr. \thevolumenumber\s
(\thevolumeyear) \startpage-\finishpage}
\immediate\write\gtoutfile{Comments: Published by Geometry and Topology Monographs at}
\immediate\write\gtoutfile{\s\s\s  http://www.maths.warwick.ac.uk/gt/GTMon\thevolumenumber/paper\thepapernumber.abs.html}
\immediate\write\gtoutfile{\noexpand\\}
\immediate\write\gtoutfile{}
\ifx\theasciiabstract\relax
\immediate\write\gtoutfile{\theabstract}\else
\immediate\write\gtoutfile{\theasciiabstract}\fi
\immediate\write\gtoutfile{}
\immediate\write\gtoutfile{\noexpand\\}
\immediate\write\gtoutfile{}
\immediate\closeout\gtoutfile}}  

\def\maketitlepage{\maketitlep\makeheadfile}
\let\makeshorttitle\maketitlepage
\let\maketitle\maketitlepage

\volumenumber{4}
\volumename{Invariants of knots and 3-manifolds (Kyoto 2001)}
\volumeyear{2002}
\papernumber{4}
\pagenumbers{43}{54}
\received{8 January 2002}
\revised{5 September 2002}
\accepted{5 September 2002}
\published{19 September 2002}

\usepackage{epsf,amssymb,amsmath,epic,eepic,pb-diagram}

\let\lbl\label

\theoremstyle{plain}

\theoremstyle{definition}

\theoremstyle{remark}

\newcommand{\psdraw}[2]{\begin{array}{c} \hspace{-1.3mm}
	\raisebox{-4pt}{\epsfxsize #2\epsfbox{#1.eps}}
	\hspace{-1.9mm}\end{array}}

\newlength{\standardunitlength}
\newcommand{\eepic}
	{\setlength{\unitlength}{0.03\standardunitlength}
	\begin{array}{c}  \hspace{-1.7mm}
         	\raisebox{-8pt}{{\renewcommand{\dashlinestretch}{30}
\begin{picture}(1165,1239)(0,-10)
\path(12,12)(12,1212)
\path(42.000,1092.000)(12.000,1212.000)(-18.000,1092.000)
\path(12,912)(312,912)
\path(192.000,882.000)(312.000,912.000)(192.000,942.000)
\path(12,612)(312,612)
\path(192.000,582.000)(312.000,612.000)(192.000,642.000)
\path(12,312)(312,312)
\path(192.000,282.000)(312.000,312.000)(192.000,342.000)
\put(612,537){\makebox(0,0)[lb]{$n$ legs}}
\end{picture}
}}\hspace{-1.9mm}\end{array}
}

\newcommand{\tr}{\operatorname{tr}}

\def\BZ{\mathbb Z}
\def\BQ{\mathbb Q}

\def\BC{\mathbb C}

\def\A{\mathcal A}
\def\B{\mathcal B}

\def\D{\Delta}

\def\La{\Lambda}
\def\l{\lambda}

\def\S{\Sigma}

\def\ihs{integral homology 3-sphere}
\def\qhs{rational homology 3-sphere}

\def\fti{finite type invariant}

\def\o1o{\underset{1}\ast}
\def\x1y{\underset{x1y}\ast}
\def\y1x{\underset{y1x}\ast}

\def\w{\omega}

\def\AS{\mathrm{AS}}
\def\IHX{\mathrm{IHX}}

\def\Wh{\mathrm{Wh}}

\def\Zrat{Z^{\mathrm{rat}}}

\def\Res{\mathrm{Res}}
\def\s{\sigma}
\def\longto{\longrightarrow}
\def\Lloc{\Lambda_{\mathrm{loc}}}
\def\hair{\mathrm{Hair}}
\def\Wh{\mathrm{Wh}}
\def\fg{\mathfrak{g}}
\def\alg{\mathrm{alg}}

\def\strutb#1#2#3{\overset{#1}{\underset{#2}{ 
\begin{array}{c} \vspace{0.0cm}
\uparrow 
\vspace{-0.25cm} \\        
| \vspace{-0.45cm} \\      
\bullet \vspace{0.00cm}   
\end{array} }}\! #3}

\begin{document}

\title{Periodicity of Goussarov-Vassiliev knot invariants}
\author{Stavros Garoufalidis}
\address{Department of Mathematics, University of Warwick\\Coventry, 
CV4 7AL, UK}
\email{stavros@maths.warwick.ac.uk}
\url{http://www.math.gatech.edu/\char'176stavros}
\begin{abstract}
The paper is a survey of known periodicity properties 
of finite type invariants of knots, and their applications.
\end{abstract}


\maketitle

\cl{\small\em Dedicated to the memory of M. Goussarov}


\section{Introduction}
\lbl{sec.intro}

\subsection{History}
\lbl{sub.history}

According to the story communicated to us by M. Polyak and O. Viro,
in the fall of 1988 M. Goussarov gave two talks about some invariants of knots
that behave in a certain sense polynomially with respect to modifications
of a knot. Goussarov's talk was neglected, and a few years later with 
entirely different motivation V. Vassiliev axiomatically introduced the
same invariants of knots. These Goussarov-Vassiliev invariants have come
to be known as {\em finite type invariants}. 

This expository paper is concerned with periodicity properties of the 
finite type invariants of knots. Our aim is to
focus on the major ideas and open problems; we will avoid proofs which the 
reader may find in the existing literature. 

\subsection{What is periodicity for the colored Jones function?}
\lbl{sub.per}

In order to motivate and explain these properties, we need to reverse our 
steps a few years earlier. In 1985 Jones discovered his famous knot 
polynomial, \cite{J}. 
To a knot Jones associated a polynomial, i.e., an element of 
$\BZ[q,q^{-1}]$ (here and below we normalize our invariants to equal to $1$
for the unknot). In one of its views, the Jones polynomial comes from a 
2-dimensional irreducible representation of the Lie algebra $\mathfrak{sl}_2$.
By considering for each natural number $d$ the $(d+1)$-dimensional
irreducible representation $V_d$ of $\mathfrak{sl}_2$, one obtains a sequence
of polynomials associated to a given knot. This sequence is best organized in 
the {\em colored Jones function} $J(K)$ of a knot $K$, which is a 2-parameter 
formal power series 
$$
J(K)(h,\l)=\sum_{n,m=0}^{\infty} a_{n,m}(K) h^n \l^m
$$
(see \cite{BG})
with remarkable periodicity properties. By its definition, if $\l=d$ is a 
natural number, $J(K)(h,d)$ coincides up to a change of variables with the 
Jones polynomial using $V_d$, thus for each fixed $d$ we have:
$$
J(K)(h,d) \in \BZ[e^h, e^{-h}].
$$
We can think of this as a {\em periodicity property} (that is, a  
{\em recursion relation}) for the coefficients $a_{n,m}$ of $J$. 
By this we mean the following: we call a sequence $(b_n)$ {\em periodic}
if the generating
function $\sum_n b_n t^n$ is a rational function of $e^t$.
Thus, for each fixed $d$, the sequence $(\sum_m a_{n,m} d^m )$ is periodic.
From our point of view, this is an obvious periodicity property of $J$.

We now discuss a {\em hidden periodicity} property of $J$. Each 
coefficient $a_{n,m}$ is a finite type invariant of degree $n$ and vanishes 
if $m > n$, \cite{BG}. Thus, we can rewrite the colored Jones function as a 
sum of subdiagonal terms $Q_{J,k}$
$$
J(K)(h,\l)=\sum_{k=0}^\infty h^k Q_{J,k}(K)(h \l)
$$
where
$$
Q_{J,k}(K)(h)=\sum_{m=0}^\infty a_{k+m,m}(K) h^m.
$$ 
It was conjectured by Melvin-Morton and Rozansky, and proven in \cite{BG}, 
that $Q_{J,0}$ equals (up to a logarithm and a change of variables) to the 
{\em Alexander polynomial} $\D$ of the knot. More precisely, we have
$$
Q_{J,0}(K)(h)=-\frac{1}{2} \log \D(K)(e^h).
$$
This translates to a hidden periodicity of the coefficients $a_{n,n}$ of $J$.

Rozansky, in his study of the colored Jones function conjectured that
for every $n \geq 1$, $Q_{J,n}$ is (up to a change of variable $t=e^h$) 
a rational function of $t$, whose denominator is a $(2n+1)$st power of the 
Alexander polynomial. Using a variety of quantum field theory techniques and 
an appropriate expansion of the $R$-martix, Rozansky proved this Rationality 
Conjecture for the colored Jones fucntion, see \cite{Ro1}.

\subsection{What is periodicity for the Kontsevich integral?}
\lbl{sub.per2}

It follows from work of Drinfeld \cite{Dr}
that the colored Jones function is the image of a 
universal finite type invariant, the so-called {\em Kontsevich integral}
$$
Z: \mathrm{Knots} \longto \A(\star)
$$ where $\A(\star)$ is a vector space over $\BQ$ 
spanned by graphs, modulo the $\AS$ and $\IHX$ relations, \cite{B}. 
The graphs in question are uni-trivalent graphs (with vertex orientations)
and have two notions of complexity: the degree (i.e., the number of vertices) 
and the number of trivalent vertices. The Kontsevich integral has a 
subdiagonal 
expansion (in terms of graphs with a fixed number of trivalent vertices).
Each term of that expansion is a series of uni-trivalent graphs with
fixed number of trivalent vertices and an arbitrary number of univalent ones.
Rozansky conjectured that each term of the expansion of the Kontsevich
integral should be given by a series of trivalent graphs with rational 
functions attached to their edges. This is often called the {\em Rationality
Conjecture} for the Kontsevich integral, which we now describe.

A weak form of the Rationality Conjecture was proven by Kricker 
\cite{Kr1}. A stong form followed by joint work with Kricker \cite{GK1}, where 
a {\em rational form} $\Zrat$ of the Kontsevich integral $Z$ was constructed.

The construction of the $\Zrat$ invariant is rather involved. A key
ingredient is the {\em surgery view of knots}, that is the presentation of
knots by surgery on framed links in a solid torus, modulo some kirby-type
relations, \cite{GK3}. Another key ingredient is that the values of the
$\Zrat$ invariant are trivalent graphs {\em with beads}, that is rational
functions attached on their edges. More precisely, 
$$
\Zrat: \, \mathrm{Knots} \longto \B(\La\to\BZ) \times \A(\Lloc)
$$
where 
\begin{itemize}
\item
$\B(\La\to\BZ)$ is a quotient of the set of Hermitian matrices over 
$\La=\BZ[t,t^{-1}]$ which are invertible over $\BZ$, modulo the equivalence
$A \sim B$ iff $A \oplus D=P(B\oplus E)P^\star$ for diagonal matrices $D,E$
with monomials in $t$ on the diagonal and for $P$ invertible over $\La$
\item 
$\A(\Lloc)$ is a vector space spanned by trivalent graphs with rational 
functions (that is elements of $\Lloc=\{f/g \, | f,g \in \La, \, g(1)=1 \}$), 
attached on their edges, modulo some relations explained in 
\cite{GK1}. 
\end{itemize}

There is a {\em hair} map 
$$
\hair: \B(\La\to\BZ) \times \A(\Lloc) \longto \A(\star)
$$
defined by
$$
\hair(A,s)= \exp\left(-\frac{1}{2} \tr\log(A)(e^h)|_{h^n \to
w_n} \right) \, \hair(s)
$$
for $A \in \B(\La\to\BZ)$ and $s \in \A(\Lloc)$, where $w_n$ is the wheel 
with $n$ legs, multiplication is given by the disjoint union and $\hair(s)$ 
replaces $t$ in each bead of $s$ by an exponential of legs as follows:
$$
\strutb{}{}{t} \to \sum_{n=0}^\infty \frac{1}{n!} \eepic 
$$
In other words, the $\hair$ map sends the matrix part to wheels and replaces 
the variable $t$ of a bead in terms of the exponential of legs. The content 
of \cite[Theorem 1.3]{GK1} is the following commutative diagram
$$
\divide\dgARROWLENGTH by1
\begin{diagram}
\node{\mathrm{Knot}}
\arrow{e,t}{\Zrat}\arrow{se,t}{Z}
\node{\B(\La\to\BZ) \times \A(\Lloc)}
\arrow{s,r}{\hair} \\
\node[2]{\A(\star)}
\end{diagram}
$$
which is in a strong form the Rationality Conjecture for the Kontsevich
integral (a weak form of the conjecture only states that the image of $Z$
is in the image of the $\hair$ map). It was recently shown by Patureau-Mirand
that the $\hair$ map is not 1-1, \cite{PM}, thus the strong form of the 
Rationality Conjecture is potentially stronger than the weak form.

\section{Properties and applications of the $\Zrat$ invariant}
\lbl{sec.app}

Although the $\Zrat$ invariant is a rather complicated object, certain parts
of it can be explained using classical topology.

\subsection{The matrix part of $\Zrat$ and the Blanchfield pairing}

In \cite{GK1} it was shown that the matrix part of $\Zrat(M,K)$ determines
the {\em Blanchfield pairing} of $(M,K)$, that is the intersection form
$$
H_1(\widetilde{M-K},\BZ) \times H_1(\widetilde{M-K},\BZ) \longto \Lloc/\La.
$$ 
The converse is also true, and will be postponed to a future publication.

\subsection{The loop move, \fti s and $0$-equivalence}
\lbl{sub.loop}

One of the beautiful ideas originating in the work of Goussarov (and Habiro)
is that a {\em geometric move} on a set $S$ leads, upon iteration, to the dual 
notions of {\em \fti s} (i.e., certain numerical invariants $f:S\to A$ with
values in an abelian group $A$ (eg. $A=\BZ,\BQ$)) and 
$n$-{\em equivalence} (i.e., certain quotients of the set $S$ that are often
abelian groups).

A main example of this idea is to consider the set of knots in $S^3$
and the move to be an $I$-{\em modification} (this is a special case of 
surgery on a clasper, and has the same result as a crossing change of a knot).
This leads to the theory of \fti s of knots in $S^3$. $0$-equivalence (that 
is knots, modulo $I$-modifications) is trivial since every knot can be 
unknotted by a sequence of $I$-modifications. 

Another example of this idea is to consider the set of knots $K$ in \ihs s $M$,
and the {\em loop move} introduced at \cite{GR}. The loop move replaces a 
pair $(M,K)$ by the result of surgery on a {\em clasper} $G \subset M-K$ 
whose leaves have linking number zero with $K$.  

The graph-part of $\Zrat$ takes values in a graded vector space (where the
degree of a trivalent graph is the number of trivalent vertices), thus
we can talk about the degree $2n$ part $\Zrat_{2n}$ of $\Zrat$. Although 
$\Zrat_{2n}$ is not a  \fti\  (since it determines a power series of finite 
type 
invariants under the $\hair$ map), it is a \fti\ of type $2n$ with
respect to the loop move, as was shown in \cite{GK1}.

In \cite{GR} it was shown that for two pairs $(M,K)$ and $(M',K')$ 
the following are equivalent:

\begin{itemize}
\item
They can be obtained one from the other via a sequence of loop moves
(this is often called $0$-{\em equivalence} with respect to the loop move).  
\item
They are $S$-{\em equivalent}, i.e., their Seifert matrices are $S$-equivalent.
\item
They have isometric Blanchfield pairings.
\end{itemize}

In conjunction with the above, it follows that the matrix part of the 
$\Zrat$ invariant precisely describes $0$-equivalence with respect to the loop 
move. 

We will now sample some applications of the $\Zrat$ invariant.


\subsection{The Casson-Walker invariant of cyclic branched coverings of a knot}
\lbl{sub.cyclic}

This application involves the computation of the Casson-Walker invariant
of cyclic branched covers of a knot. 
Given a knot $K$ in an \ihs\ $M$ and a positive integer $p$, one can construct
a closed 3-manifold, the $p$-fold branched covering $\S^p_{(M,K)}$ of $K$
in $M$. The sequence of 3-manifolds $\S^p_{(M,K)}$ is closely related to the
{\em signature function} $\s(M,K): S^1 \to \BZ$ of a knot, which
properly normalized is a key concordance invariant. The first nontrivial
finite type invariant of 3-manifolds is the Casson-Walker invariant $\l$.
Suppose that $\S^p_{(M,K)}$ is a \qhs\ (or equivalently, that the Alexander
polynomial of $(M,K)$ has no complex $p$th roots of unity). 
In \cite[Corollary 1.4]{GK2} we showed that for all $(M,K)$ and $p$ as above, 
we have
$$
\l (\S^p_{(M,K)})= \frac{1}{3} \Res_p \Zrat_2(M,K)
+ \frac{1}{8} \s_p(M,K)
$$
where $\s_p(M,K)=\sum_{\omega^p=1} \s(M,K)(\omega)$ is the {\em total} 
$p$-{\em signature} and 
$$
\Res_p \left( \psdraw{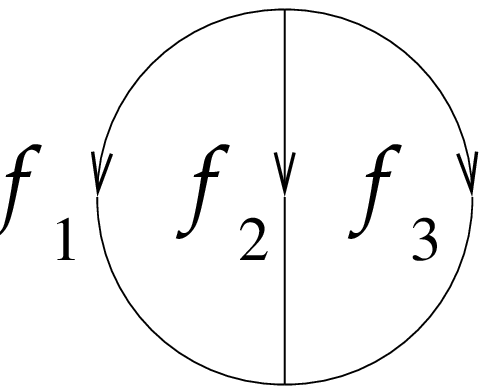}{0.6in} \right)= \frac{1}{p} \sum 
f_1(\w_1) f_2(\w_2) f_3(\w_3)
$$
where the sum is over all triples $(\w_1,\w_2,\w_3)$ that satisfy
$\w_1^p=\w_2^p=\w_3^p=\w_1 \w_2 \w_3=1$. In particular, the growth rate 
of the Casson invariant of cyclic branched covers is given by:
$$
\lim_{p \to \infty} \frac{\l (\S^p_{(M,K)})}{p} = \frac{1}{3}
\int_{S^1 \times S^1} Q(M,K)(s) d\mu(s) +
\frac{1}{8} \int \s_s(M,K) d\mu(s)
$$
where $d\mu$ is the Haar measure. The above formulas are part of the
computation of the full LMO invariant of cyclic branched covers in terms
of the signature function and residues of the $\Zrat$ invariant, 
\cite[Theorem 1]{GK2}.

\subsection{On knots with trivial Alexander polynomial}
\lbl{sub.good}

Topological surgery in dimension $4$ is a list of surgery problems that
sometimes one can be reduced to another. A list of {\em atomic surgery 
problems} (that is problems that every surgery problem reduces to, in dimension
$4$) were compiled by Casson and Freedman. Such a list is not unique but
various versions of it involve the {\em slicing of boundary links} whose free 
covers are acyclic. In the case of knots, this condition is equivalent to
the triviality of the Alexander polynomial. Freedman showed that
Alexander polynomial $1$ knots are indeed topologically slice, \cite{F2}.

However, for over 15 years there was a confusion: it was thought that
every Alexander polynomial $1$ knot bounds a Seifert surface whose Seifert
form has minimal rank (that is, rank equal to the genus of the surface),
\cite{F1}. 
This turnt out to be false, \cite{GT1}. The required invariant turnt out to 
be the 2-loop part $\Zrat_2$ of the rational invariant $\Zrat$. No classical
knot invariant (such as signature,
Blanchfield pairing, Casson-Gordon, and Cochran-Orr-Teichner $L^2$ signature
invariants) would have worked, as all these vanish on knots with trivial
Alexander polynomial. Thus, the $\Zrat_2$ invariant is an obstruction
to a knot bounding a Seifert surface of a prescribed type.

In subsequent work with P. Teichner, we introduce a
decreasing filtration on the set of Alexander polynomial $1$ knots 
\cite{GT2}, which is strictly decreasing every other step, and which in
degree $1$ equals to the set of knots which bound minimal rank Seifert 
surfaces. The
{\em Hyperbolic Volume Conjecture} of knots (in its Simplicial Volume 
formulation) implies that the intersection of this filtration is the unknot,
\cite{MM}.

Whether this descreasing filtration on good knots is related to a tower of
smooth slicing obstructions of good knots is an unknown problem.

\subsection{Lifting from the Lie algebra to the Lie group}
\lbl{sub.lie}

Given a Lie algebra $\fg$ with invariant inner product, there is a map
(often called a {\em weight system})
$$
W_{\fg}^h: \A(\star) \longto S(\fg)^{\fg}[[h]]
$$
where $S(\fg)$ is the symmetric algebra of $\fg$. The image of the Kontsevich
integral under this map coincides with the colored Jones function. This
weight system replaces unitrivalent graphs by Lie invariant tensors.

Similarly, given a compact connected Lie group $G$ with (complexified)
Lie algebra $\fg$, there is a map
$$
W_G^h: \A(\Lloc) \longto C_{\alg}(G)^G[[h]]
$$
where $C_{\alg}(G)$ is the algebra of {\em almost invariant functions}
$f: G \to \BC$ (that is those continuous functions that generate a finite 
dimensional
subspace of the algebra of continuous functions $C(G)$ under the action of $G$
on itself by conjugation). $C_{\alg}(G)^G$ is a finitely generated algebra,
as follows from the Peter-Weyl theorem. For a discussion, see \cite{Ga2}.

Using this, it is fairly easy to show that the Rationality Conjecture
for the Kontsevich integral implies the Rationality Conjecture for the 
colored Jones function, see \cite[Theorem 2]{Ga2} and also \cite{Ro2}. 

\subsection{Sample calculations}
\lbl{sub.sample}

Can we compute the $\Zrat$ of any knot? This is a hard question, since such
a computation would in particular compute the Kontsevich integral of that knot.
Computations can be done for torus knots, generalizing work of Bar-Natan
and Lawrence, \cite{BL}. On the other hand, the author does not know how to
compute the $\Zrat$ invariant of the figure eight knot.

One can ask for less, however. Rozansky has written a computer program that
can compute the $\Zrat_2$ invariant of knots which are presented as closures
of braids, \cite{Ro2}

A different algorithm for computing $\Zrat_2$ for knots has been given
in \cite{GR}, where knots are presented via surgery on claspers.
 For more computations, see also \cite{GT1} and \cite{Kr2}. 
We mention a sample 
computation here for an untwisted Whitehead double $\Wh(K)$ of a knot $K$, 
taken from \cite[Corollary 1.1]{Ga1}:
$$
\Zrat_2(\Wh(K))=a(K) \cdot \psdraw{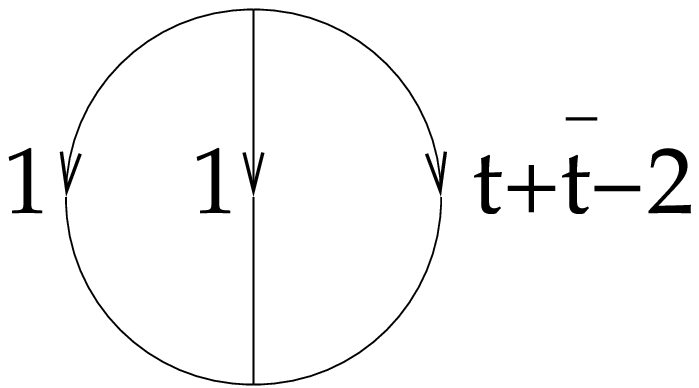}{1in}
$$
where $a(K)=\frac{1}{4}\frac{d^2}{dh^2}\D(K)(e^h)_{h=0} \in \BZ$.
This implies that the Kontsevich integral of $\Wh(K)$ is given by
$$
Z(\Wh(K))=\exp\left(2 \sum_{n=1}^\infty \frac{a(K)}{(2n)!} 
\psdraw{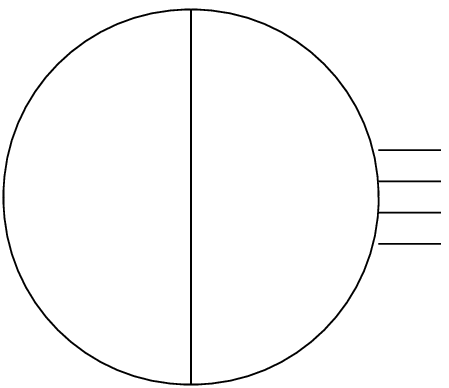}{0.5in} \}_{2n \, \mathrm{legs}} +...\right)
$$
where the dots indicate connected unitrivalent graphs with at least three
loops and an arbitrary number of legs. This illustrates the periodicity
of the 2-loop part of the Kontsevich integral.

\subsection{Integrality properties of the colored Jones function}
\lbl{sub.integrality}

This is a topic complementary to periodicity, that we briefly discuss.
Recall that $J(h,d) \in \BZ[e^h,e^{-h}]$, which is an {\em integrality
property} of coefficients $a_{n,m}$ of $J$, in the sense that certain
linear combinations of rational numbers are actually integers.
Rozansky conjectured and further proved in \cite{Ro1} that each of the
subdiagonal terms $Q_{J,k}$ ($k \geq 1$) of the colored Jones polynomial  
can be written in the form
$$
Q_{J,k}(h)=\frac{P_{J,k}(e^h)}{\D^{2n+1}(e^h)}
$$
for polynomials $P_{J,k}(t) \in \BZ[t,t^{-1}]$. This is a hidden integrality
property of the coefficients of $J$, and can be used, as Rozanky showed,
to show that the image of the LMO invariant under the $\mathfrak{sl}_2$
weight system coincides with the $p$-adic expansion of the {\em 
Reshetikhin-Turaev} invariant of manifolds obtained by nonzero surgery on a 
knot, \cite{Ro3}. 

\section{What next?}
\lbl{sec.what}

\subsection{Boundary links}
\lbl{sub.blinks}
The Rationality Conjecture (RC) for the Kontsevich integral of knots has been 
extended to links in two different directions. 
In one extension, one considers boundary links. If we are to have a RC
for the Kontsevich integral of links, we need to restrict to links whose
Kontsevich integral has vanishing tree part. This is equivalent to the
vanishing of all Milnor invariants, as was shown by Habegger-Masbaum, 
\cite{HM}. The class of links with vanishing Milnor invariants contains
(and perhaps coincides with) the class of {\em homology boundary links},
\cite{L}. For simplicity, we can restrict our attention to boundary
links.

An extra point to keep in mind is that the Kontsevich integral of a link
takes value in uni-trivalent graphs whose legs are colored by the components
of the link. Further, legs with different colored do not commute (i.e.,
the order by which they are attached on an edge is important). In this
context, the RC for the Kontsevich integral of boundary links is a statement
about {\em rational functions in noncommuting variables}. Luckily, there
is well-developed algebra to deal with this, that comes under the name
of {\em noncommutative localization} (i.e., Cohn localization) of the
group-ring of the free group. Adapting this, it was possible to define a 
rational noncommutative invariant $\Zrat$ of boundary links that 
determines the Kontsevich integral of the boundary link.

\subsection{Algebraically connected links}
\lbl{sub.algc}

In a direction perpendicular to that of boundary links, Rozansky considers
links $L$ with nonvanishing (multivariable) Alexander polynomial.  
The simplest example of such a link is the Hopf link.

Rozansky's RC is a {\em relative} conjecture, formulated for the colored Jones
function. Whether it can be lifted to a relative RC for the Kontsevich integral
of alg. connected links remains to be seen.

\subsection{More periodicity properties?}
\lbl{sub.more}

The best that one could hope for is that after a change of variables, the
two parameter series $J(h,\l)$ is given by a single rational function.
This is probably false, but can we prove it?

\subsection{Acknowledgements}
The author presented various parts of this paper in several occasions
including the DennisFest (Stony Brook, June 2001), the International Topology
Conference in Athens Georgia (May 2001), the LevineFest (Tel-Aviv, 2001),
the Goussarov Day in Kyoto (September 2001), a colloqium at Liverpool 
(December 2001) and a series of lectures at Paris 7 (December 2001). The
author wishes to thank M. Farber, G. Matic, G. Masbaum, H. Morton, 
D. Sullivan, P. Vogel, S. Weinberger and especially H. Murakami and 
T. Ohtsuki for their hospitality, and the persistent audience for their 
remarks. 

This research was partially supported by NSF and BSF.

\Addresses\recd

\end{document}

\thanks{

\date{
This edition: September 5, 2002 \hspace{0.3cm} 
First edition: January 8, 2002.}